\documentclass{article}
\usepackage{amsmath,amssymb}
\usepackage{theorem}
\usepackage{color}

\theorembodyfont{\itshape}
\theoremstyle{plain}
\newtheorem{theorem}{Theorem}[section]
\newtheorem{lemma}{Lemma}[section]
\newtheorem{corollary}{Corollary}[section]
\newtheorem{proposition}{Proposition}[section]
\theorembodyfont{\rmfamily}

\newtheorem{construction}{Construction}[section]

\newtheorem{remark}{Remark}[section]

\newcommand{\Bd}{\mbox{{\rm Bd}}\,}

\newcommand{\Ind}{\mbox{{\rm Ind}}\,}
\newcommand{\ind}{\mbox{{\rm ind}}\,}

\begin{document}

\title{\bf On a modification of finite-dimensional Niemytzki spaces.}

\author{Vitalij A.~Chatyrko}


\maketitle

\begin{abstract} In this paper we extend the construction of the Niemytzki plane to dimension $n \geq 3$. Further, we consider a poset of topologies on the closed $n$-dimensional Euclidean half-space
similar to one  from \cite{AAK}  which is related to the Niemytzki plane topology. Then we explore properties of topologies from the poset.

\end{abstract} 

\medskip
{\it Keywords and Phrases: a finite-dimensional Niemytzki space; a Lindel\"of space; a perfect space: a $\sigma$-compact space}  

\smallskip
{\it 2000 AMS (MOS) Subj. Class.:} Primary 54A10, 54D15   
\medskip
\baselineskip=18pt

\section{Introduction}
The Niemytzki plane (cf. \cite[Example 82]{SS}) is a classical example of 
a topological space (like the square of the Sorgenfrey line (cf. \cite[Example 84]{SS})) which is Tychonoff but not normal.  Besides that the Niemytzki plane is a separable, first-countable, perfect,  realcompact, Cech-complete space which is neither countably paracompact nor weakly paracompact (cf. \cite{E}).
Recently it was proved that the Niemytzki plane is even $\kappa$-metrizable {\cite{BKP}}.

In this paper we extend the construction of the Niemytzki plane to dimension $n \geq 3$. Further, we consider a poset of topologies on the closed $n$-dimensional Euclidean half-space
similar to one  from \cite{AAK}  which is related to the Niemytzki plane topology. Then we explore properties of topologies from the poset.

For standard notions we refer to \cite{E}.

\section{Finite-dimensional Niemytzki spaces and their properties}
We generalize to dimension $n \geq 3$ the construction of Niemytzki plane (cf. \cite[Example 82]{SS}).

\begin{construction} Consider subsets $P_n = \{\overline{x}=(x_1, \dots, x_n): x_i \in \mathbb R, x_n > 0\}$ and $L_n = \{\overline{x}= (x_1, \dots, x_n): x_i \in \mathbb R, x_n = 0\}$ of $\mathbb R^n$. 

We generate a topology $\tau_N$ on $X_n = P_n \cup L_n$ as follows.

If $\overline{a} \in P_n$ then a local base of $\tau_N$ at $\overline{a}$ consists of sets 
$B(\overline{a}, \epsilon) = \{\overline{x} \in R^n : |\overline{x}-\overline{a}| < \epsilon\}$, where
$|\overline{x}-\overline{a}| = \sqrt{\Sigma_{i=1}^n (x_i -a_i)^2}$ and $0 < \epsilon < a_n$.

If $\overline{a} \in L_n$ then a local base of $\tau_N$ at $\overline{a}$ consists of sets 
$\tilde{B}(\overline{a}, \epsilon) = \{\overline{a}\} \cup B(\overline{a(\epsilon)}, \epsilon)$, where $\overline{a(\epsilon)} = (a_1, \dots, a_{n-1}, \epsilon)$
 and $0 < \epsilon$.
 \end{construction}
 
Let us list some obvious properties of the finite-dimensional Niemytzki spaces.

Note that the space $(X_n, \tau_N)$ is first-countable, the restriction of the topology $\tau_N$ onto $P_n$ (resp. $L_n$) coincides with the Euclidean topology on $P_n$ (resp. with the discrete topology on $L_n$ ), the subset $P_n$ (resp. $L_n$) is open and dense (resp. closed and nowhere dense) in the space $(X_n, \tau_N)$. In particular, the space $(X_n, \tau_N)$ is separable and perfect. 

Denote the Euclidean topology on $X_n$ by $\tau_E$. It is evident  that $\tau_E \subseteq \tau_N$.
So the space $(X_n, \tau_N)$ is completely Hausdorff.

It is easy to see that for any $\overline{a} \in L_n$ the boundary 
$\Bd_N \tilde{B}(\overline{a}, \epsilon)$ of
the set $\tilde{B}(\overline{a}, \epsilon)$ in the space $(X_n, \tau_N)$ is equal to 
$\Bd_E B(\overline{a(\epsilon)}, \epsilon) \setminus \{\overline{a} \}$, where 
$\Bd_E B(\overline{a(\epsilon)}, \epsilon)$ is the boundary of the set  
$B(\overline{a(\epsilon)}, \epsilon)$ in the space $(X_n, \tau_E)$.

Moreover, $B(\overline{a(\epsilon)}, \epsilon)$ is the disjoint union of the sets 
$\Bd_N \tilde{B}(\overline{a}, t\epsilon)$, where $0<t<1$. Note that for each $\overline{x} \in B(\overline{a(\epsilon)}, \epsilon)$ the point $\overline{x}$ belongs to the only set
$\Bd_N \tilde{B}(\overline{a}, t(\overline{x})\epsilon)$ with 
$t(\overline{x}) = (\Sigma_{i=1}^{n-1} (x_i - a_i)^2 + x_n^2)/(2 \epsilon x_n).$

\begin{proposition} The space $(X_n, \tau_N)$ is Tychonoff.
\end{proposition}

Proof. Since $\tau_E \subseteq \tau_N$, it is enough to show that for any $\overline{a} \in L_n$ and any $\tilde{B}(\overline{a}, \epsilon)$  there exists a continuous function $f: (X_n, \tau_N) \to [0,1]$ such that 
$f(\overline{a}) = 0$ and $f|_{\Bd_N \tilde{B}(\overline{a}, \epsilon)} = 1$. In fact, 
set $f(\overline{x}) = 1$ for any $\overline{x} \in X_n \setminus \tilde{B}(\overline{a}, \epsilon)$,
$f(\overline{a}) = 0$  and $f(\overline{x}) = t(\overline{x})$ for any  $\overline{x} \in B(\overline{a(\epsilon)}, \epsilon)$. $\Box$

It is easy to see that the set $X'_{m,n} = \{\overline{x} \in X_n : x_1 = \dots x_{n-m} =0\}$, where $2 \leq m < n$, is a closed
subset of $(X_n, \tau_N)$, and the subspace $X'_{m,n}$ of $(X_n, \tau_N)$ is homeomorphic to the space 
$(X_{m}, \tau_N)$. Since the space $(X_2, \tau_N)$ is the Niemytzki plane, the space 
$(X_n, \tau_N)$ is in particular neither normal, countably paracompact nor weakly paracompact.

\section{Topologies between the Euclidean and finite-dimensional Niemytzki}

In \cite{AAK} Abuzaid, Alqahtani and Kalantan suggested by the use of technique from \cite{H} a poset $\mathcal T$ of topologies on the set $X_2$ such that the minimimal topology is $\tau_E$ and the maximal topology is $\tau_N$.
We will extend the construction to the sets $X_n, n \geq 3$.

\begin{construction} Let $A$ be a subset of $L_n$. We generate a topology $\tau(A)$ on $X_n = P_n \cup L_n$ as follows.

If $\overline{a} \in P_n$ then a local base of $\tau(A)$ at $\overline{a}$ consists of sets 
$B(\overline{a}, \epsilon)$, where $0 < \epsilon < a_n$.

If $\overline{a} \in A$ then a local base of $\tau(A)$ at $\overline{a}$ consists of sets 
$B(\overline{a}, \epsilon) \cap X_n$, where $0 < \epsilon$.

If $\overline{a} \in L_n \setminus A$ then a local base of $\tau(A)$ at $\overline{a}$ consists of sets 
$\tilde{B}(\overline{a}, \epsilon)$, where $0 < \epsilon$.

\end{construction}

It is evident  that $\tau_E = \tau(L_n) \subseteq  \tau(A) \subseteq \tau(\emptyset) = \tau_N$. Note also that for any  $A, B \subseteq L_n$ we have  $A \subseteq B$ iff $\tau(A) \supseteq \tau(B)$.

Note that the space $(X_n, \tau(A))$ is first-countable and separable, the restriction of the topology $\tau(A)$ onto $P_n \cup A$ (resp. $L_n \setminus A$) coincides with the Euclidean topology on $P_n \cup A$ (resp. with the discrete topology on $L_n \setminus A$ ), the subset $P_n$ (resp. $L_n$) is open and dense (resp. closed, even a zero set,  and nowhere dense) in the space $(X_n, \tau(A))$.

Similarly to  $(X_n, \tau_N)$, one can prove that the space $(X_n, \tau(A))$ is Tychonoff.

\begin{proposition}(for $n=2$ see \cite[Theorem 2.3]{AAK}) The following are equivalent.
\begin{itemize}
\item[(i)] The space $(X_n, \tau(A))$ is hereditarily Lindel\"of. 
\item[(ii)] $|L_n \setminus A| \leq \aleph_0$.
\item[(iii)] The space $(X_n, \tau(A))$ is second-countable.
\item[(iv)] The space $(X_n, \tau(A))$ is metrizable.
\end{itemize}
\end{proposition}

Proof.  $(i) => (ii)$.  Since the space $(L_n \setminus A, \tau(A)|_{L_n \setminus A})$ is discrete, $|L_n \setminus A| \leq \aleph_0$.
 $(ii) => (iii)$. Let $L_n \setminus A = \{\overline{b_1}, \overline{b_2}, \dots\}$,
$\mathcal B$ be a countable base for the space  $(X_n, \tau_E)$ and $\mathcal B_i$ be a countable 
local base of the space $(X_n, \tau(A))$ at the point $\overline{b_i}, i = 1, 2, \dots$  
Then the family $\mathcal B \cup \cup_{i=1}^\infty \mathcal B_i$ is a  countable base  for
the space $(X_n, \tau(A))$.
$(iv) => (i)$. Since the space $(X_n, \tau(A))$ is separable, it is second-countable and hence hereditarily Lindel\"of. $\Box$

Let $\overline{a} \in L_n \setminus A$. It is easy to see that any sequence of points 
$\{\overline{x_i}\}_{i=1}^\infty$ in $\Bd_E B(\overline{a(\epsilon)}, \epsilon) \setminus \{\overline{a} \}$ converging to $\overline{a}$ in the Euclidean topology $\tau_E$ 
is discrete in the space $(X_n, \tau(A))$. This implies the following proposition.

\begin{proposition}(for $n=2$ see \cite[Theorem 2.6]{AAK}) The space $(X_n, \tau(A))$ is locally compact iff $A = L_n$ i. e. $\tau(A) =\tau_E$.
\end{proposition}

\section{Some other properties of the spaces $(X_n, \tau(A))$}

The following proposition is evident.

\begin{proposition}\label{prop_1} Let $A$ be any subset of $L_n$. Then 
$(X_n, \tau(A))$ is perfect (resp. Lindel\"of or $\sigma$-compact) iff
$(L_n, \tau(A)|_{L_n})$ is the same.
\end{proposition}

Let us observe that the topology $\tau(A)|_{L_n}$
can be considered as a modification of the Euclidean topology on the set $L_n$ in the sense of Bing \cite{B}
and Hanner \cite{Han}, see \cite[Example 5.1.22]{E} for the general construction.

Note also that the space $(L_n, (\tau_E)|_{L_n})$ is homeomorphic to the Euclidean space $\mathbb R^{n-1}$, $(\tau_E)|_{L_n} \subseteq \tau(A)|_{L_n}$, $(\tau_E)|_A = \tau(A)|_A$, the set
$A$ is closed in the space  $(L_n, \tau(A)|_{L_n})$ and the subspace $L_n \setminus A$ of $(L_n, \tau(A)|_{L_n})$ is discrete.

So  from \cite[Problem 5.5.2 (b)]{E} we get that the space $(L_n, \tau(A)|_{L_n})$ is hereditarily collectionwise normal.

By 
\cite[Problem 5.5.2 (c)]{E} we easily obtain

\begin{lemma}\label{cor_1} The space $(L_n, \tau(A)|_{L_n})$ is perfect iff  $A$  is a $G_\delta$-set in  $(L_n, (\tau_E)|_{L_n}).$
\end{lemma}

\begin{lemma}\label{cor_2} The space $(L_n, \tau(A)|_{L_n})$ is Lindel\"of iff 
$L_n \setminus A$ does not contain a closed uncountable subset of 
 $(L_n, (\tau_E)|_{L_n})$.
\end{lemma} Proof. $=>$ Let us note that any closed subset $Y$ of 
 $(L_n, (\tau_E)|_{L_n})$ such that $Y \subset L_n \setminus A$ is a closed discrete subset of $(L_n, \tau(A)|_{L_n})$.
 Hence, $|Y|\leq \aleph_0$.

$<=$ Let $\alpha$ be an open cover of the space $(L_n, \tau(A)|_{L_n})$. Since  $\tau(A)|_A = (\tau_E)|_A$, there exists a countable subfamily $\alpha_1$ of $\alpha$  such that 
$A \subseteq \cup \alpha_1$. Moreover, by the definition of $\tau(A)$ there is 
an open set $O$ of  $(L_n, \tau_E)$ such that 
$A \subseteq O \subseteq \cup \alpha_1$.
Let us note that $B= L_n \setminus O  \subseteq L_n \setminus A$ is a closed subset of 
$(L_n, \tau_E)$. By assumption $B$ is countable. So there exists a countable subfamily $\alpha_2$ of $\alpha$ such that $B \subseteq \cup \alpha_2$. 
Let us observe that the subfamily $\alpha_1 \cup \alpha_2$ of $\alpha$ is countable and it covers $L_n$. 
$\Box$

\begin{lemma}\label{cor_3} The space $(L_n, \tau(A)|_{L_n})$ is $\sigma$-compact iff 
 $A$  is a $F_\sigma$-set in  $(L_n, (\tau_E)|_{L_n})$ and $|L_n \setminus A| \leq \aleph_0$.
\end{lemma}
Proof. $=>$ Since $A$ is closed in $(L_n, \tau(A)|_{L_n})$, $A$ is $\sigma$-compact in 
$(L_n, \tau(A)|_{L_n})$ and hence in $(L_n, \tau(E)|_{L_n})$. So  $A$  is a $F_\sigma$-set in  $(L_n, (\tau_E)|_{L_n})$. Further, $L_n \setminus A$ is a $G_\delta$-set in $(L_n, (\tau_E)|_{L_n})$ which is homeomorphic to $\mathbb R^{n-1}$. If $L_n \setminus A$ is uncountable then there is an uncountable compact subset $Y$ of 
$(L_n, (\tau_E)|_{L_n})$ such that $Y \subset L_n \setminus A$. Let us note that $Y$ is a closed discrete
uncountable
subset of $(L_n, \tau(A)|_{L_n})$. Since the space $(L_n, \tau(A)|_{L_n})$ is $\sigma$-compact it is impossible. So $|L_n \setminus A| \leq \aleph_0$.

$<=$ It is trivial. $\Box$ 

Proposition \ref{prop_1} and Corollories \ref{cor_1}-\ref{cor_3} imply
\begin{theorem}\label{theo_1} Let $A$ be any subset of $L_n$. Then 
$(X_n, \tau(A))$ is perfect (resp. Lindel\"of or $\sigma$-compact) iff
 $A$  is a $G_\delta$-set in  $(L_n, (\tau_E)|_{L_n})$ (resp.
 $L_n \setminus A$ does not contain a closed uncountable subset of 
 $(L_n, (\tau_E)|_{L_n})$ or  $A$  is a $F_\sigma$-set in  $(L_n, (\tau_E)|_{L_n})$ and $|L_n \setminus A| \leq \aleph_0$
 ).
\end{theorem}

\begin{corollary} \label{cor} Let $B \subseteq L_n \setminus A$ be a closed uncountable subset of $(L_n, (\tau_E)|_{L_n})$,
for example, $B$ is homeomorphic to the Cantor set. Then the space $(X_n, \tau(A))$ is not Lindel\"of.
\end{corollary}

\begin{remark} Since any Cantor set in $\mathbb R^n$ is nowhere dense in $\mathbb R^n$, Corollary \ref{cor} (the case $n=2$) evidently disagrees with \cite[Theorem 2.8 and Theorem 2.9]{AAK}.
Indeed, the fourth sentence of the proof \cite[Theorem 2.8]{AAK} is not correct. Roughly speaking it implies that any open dense subset of the real line must coincide with the real line. 
\end{remark}

\begin{corollary} If $(X_n, \tau(A))$ is $\sigma$-compact then  $(X_n, \tau(A))$  is second-countable.
\end{corollary}

\begin{corollary} Let $M$ be any countable dense subset of $(L_n, (\tau_E)|_{L_n})$. Then
\begin{itemize}
\item[(1)] $(X_n, \tau(M))$ is neither perfect nor Lindel\"of,
\item[(2)] $(X_n, \tau(L_n \setminus M))$ is second-countable but it is not $\sigma$-compact.
\end{itemize}

\end{corollary} 

\begin{corollary} Let $A$ be any uncountable compact subset of  $(L_n, (\tau_E)|_{L_n})$. Then 
$(X_n, \tau(A))$  (as well as $(X_n, \tau(L_n \setminus A))$) is perfect but it is not Lindel\"of.
\end{corollary}

A subset $A$ of the Euclidean space $\mathbb R^n$ we will call {\it a Bernstein set} (cf. \cite[p. 24]{O} for the case $n=1$)  if both $A$ and
$\mathbb R^n \setminus A$  intersect every uncountable compact subspace $F$
 of $\mathbb R^n$. It is easy to see that if $A$ is a Bernstein  set of $\mathbb R^n$ then $\mathbb R^n \setminus A$ is 
 also a Bernstein set of $\mathbb R^n$. Moreover, the Bernstein sets are of the cardinality continuum and they do not contain uncountable compacta. It implies, in particular, that any Bernstein set in  $\mathbb R^n$  is neither a $G_\delta$-set nor  an $F_\sigma$-set in $\mathbb R^n$.

\begin{corollary} \label{Lindel} Let 
$A$ be a Bernstein set of the space $(L_n, (\tau_E)|_{L_n})$. Then the space  
$(X_n, \tau(A))$ is Lindel\"of but it is not perfect. 
\end{corollary}

\begin{remark} Similarly to \cite{M} one can prove that the topological product 
$(X_n, \tau(A)) \times (L_n \setminus A, \tau(E)|_{L_n \setminus A})$  is not normal, where
$A$ is a Bernstein set of the space $(L_n, (\tau(E))|_{L_n})$.
\end{remark}

Recall (cf. \cite[p. 65]{Cha}) that a subset $Y$ of a space $X$ is called
\begin{itemize}
\item[(1)] {\it $C^*$-embedded in $X$}  if
every bounded continuous function on $Y$ can be extended to a bounded continuous function on $X$,
\item[(2)] {\it $z$-embedded} in $X$ if every zero set of $Y$ is the trace on $Y$ of some zero set of $X$.
\end{itemize}

Let us recall that any closed subset $Y$ of a normal space $X$ is   $C^*$-embedded in $X$, and if a subset $Y$ of a space $X$  is $C^*$-embedded in $X$ then $Y$ is $z$-embedded in $X$.

\begin{lemma}\label{z-embed} Let $Y$ be a discrete subspace  of cardinality continuum of a separable space $X$.
 Then $Y$ is not $z$-embedded in $X$.
\end{lemma}

Proof. Let $X$ be a separable space,  $Y$ be its discrete subspace of cardinality  $\mathfrak{c}$, $\mathcal F$ be the family of all continuous functions on the space $X$ and $Z_X$ be the family of all zero sets on $X$. 
Since $X$ is separable,   the cardinality of $\mathcal F$ is at most $\mathfrak{c}$ and hence
the cardinality of $Z_X$ is also at most $\mathfrak{c}$. 
Let $Z_Y$ be the family of all zero sets on $Y$.  
It is easy to see that  the cardinality of $Z_Y$ is at least $2^\mathfrak{c} > \mathfrak{c}$. 
So $Y$ is not $z$-embedded in $X$.
 $\Box$

The following corollary is evident.
\begin{corollary} 
No separable  normal  space contains a closed discrete subspace of cardinality continuum.
\end{corollary}

Recall (\cite[Exercise 5.2. C (b)]{E}) that no separable countably paracompact  space contains a closed discrete subspace of cardinality continuum.

\begin{proposition} \label{not normal} Let $B \subseteq L_n \setminus A$ be a closed uncountable subset of  $(L_n, (\tau_E)|_{L_n})$.
Then the space $(X_n, \tau(A))$ is neither normal nor 
countably paracompact. 
\end{proposition}

Proof. 
Since $(\tau_E)|_{L_n} \subseteq \tau(A)|_{L_n}$, the set $B$ is a closed discrete subset of 
$(L_n, \tau(A)|_{L_n})$ (and even $(X_n, \tau(A))$) of cardinality continuum. 
Since the space $(X_n, \tau(A))$ is separable, it is neither normal nor 
countably paracompact. $\Box$

\begin{corollary} \label{cor_nor} If the space $(X_n, \tau(A))$ is normal (or 
countably paracompact) then $L_n \setminus A$ does not contain a closed uncountable subset of $(L_n, (\tau_E)|_{L_n})$.
\end{corollary}

\begin{theorem}\label{equi}  The following are equivalent.

\begin{itemize}
\item[(i)] The space $(X_n, \tau(A))$ is Lindel\"of.

\item[(ii)]  The space $(X_n, \tau(A))$ is paracompact.

\item[(iii)]  The space $(X_n, \tau(A))$ is countably paracompact.

\item[(iv)]  The space $(X_n, \tau(A))$ is normal.

\item[(v)]  The set $L_n \setminus A$ does not contain a closed uncountable subset $(L_n, (\tau_E)|_{L_n})$. 

\end{itemize}

\end{theorem}

Proof. Let us only note that  $(iii) $ ( or $(iv)$)$=> (v)$ (by Corollary \ref{cor_nor}).
$(v) => (i)$ (by Theorem \ref{theo_1}). $\Box$

The following lemma is evident.
\begin{lemma} \label{z-embed_2}Let $X$ be a space and  $Z \subseteq Y \subseteq X$. If $Z$ is $z$-embedded in $Y$ and $Y$
is $z$-embedded in $X$ then $Z$ is $z$-embedded in $X$.
\end{lemma}

\begin{proposition} 
\begin{itemize}
\item[(i)] Let $B$ be a subset of $L_n \setminus A$ of cardinality continuum. Then $B$ is not 
$z$-embedded in
$(X_n, \tau(A))$. 
\item[(ii)] Let $B \subseteq L_n \setminus A$ be a closed uncountable subset of $(L_n, (\tau_E)|_{L_n})$,
for example, $B$ is homeomorphic to the Cantor set. Then the subset $L_n$ of $(X_n, \tau(A))$ is not $z$-embedded in
$(X_n, \tau(A))$. 

\end{itemize}
\end{proposition}
Proof. (i) Observe that the subspace $(B, \tau(A)|_B)$ of the space  $(X_n, \tau(A))$ is discrete and has the cardinality continuum.
Since $(X_n, \tau(A))$ is separable, by Lemma \ref{z-embed}  the set $B$ is not $z$-embedded in
$(X_n, \tau(A))$. 

(ii) Note that $B$ is a closed discrete subset of the space $(L_n, \tau(A)|_{L_n})$ (and also $(X_n, \tau(A)$) of cardinality continuum.
Hence by (i) the set $B$ is not $z$-embedded in
$(X_n, \tau(A))$. 
Recall that  the space  $(L_n, \tau(A)|_{L_n})$ is normal. So the set $B$ is  $z$-embedded in $(L_n, \tau(A)|_{L_n})$. If we assume that the closed subset $L_n$ of $(X_n, \tau(A))$ is  $z$-embedded in
$(X_n, \tau(A))$ we will get a contradiction with Lemma \ref{z-embed_2}. $\Box$

\begin{corollary} \label{z-emb} If the subset $L_n$ of $(X_n, \tau(A))$ is $z$-embedded in
$(X_n, \tau(A))$ then the set $L_n \setminus A$ does not contain a closed uncountable subset of 
$(L_n, (\tau_E)|_{L_n})$,
\end{corollary}

\begin{theorem} The following  are equivalent.
\begin{itemize} 
\item[(i)] The space $(X_n, \tau(A))$ is normal.

\item[(ii)] The subset $L_n$ of $(X_n, \tau(A))$ is $C^*$-embedded in
$(X_n, \tau(A))$.

\item[(iii)] The subset $L_n$ of $(X_n, \tau(A))$ is $z$-embedded in
$(X_n, \tau(A))$.
\end{itemize}
\end{theorem}

Proof.  $(iii) => (i)$ (by Corollary \ref{z-emb} and Theorem \ref{equi}). $\Box$

\begin{remark} For any subset $A$ of $L_n$ we have 

$\ind  (L_n, \tau(A)|_{L_n}) = \Ind  (L_n, \tau(A)|_{L_n}) = \dim  (L_n, \tau(A)|_{L_n}) = \dim A$, and

$\ind (X_n, \tau(A)) = \Ind (X_n, \tau(A))= \dim (X_n, \tau(A)) = n$ 
for any $A$ for which the space $(X_n, \tau(A))$ is normal.
\end{remark}

\noindent(V.A. Chatyrko)\\
Department of Mathematics, Linkoping University, 581 83 Linkoping, Sweden.\\
vitalij.tjatyrko@liu.se

\end{document}